# Random chemostats with competition and different kinetics to investigate the growth of the gut microbiome

Javier López-de-la-Cruz[a,1], Felipe Rivero[b] and Carlos R. Takaessu Jr[c]

[a] Dpto. de Matemática Aplicada a las TIC, Escuela Técnica Superior de Ingenieros Informáticos, Campus de Montegancedo, Universidad Politécnica de Madrid, 28660 Boadilla del Monte, Madrid, Spain.

[b] Dpto. de Matemática Aplicada a las TIC, Escuela Técnica Superior de Ingeniería y Sistemas de Telecomunicación, Campus Sur, Universidad Politécnica de Madrid, 28031, Madrid, Spain.

[c] Dpto. de Matemática, Instituto de Ciências Matemáticas e de Computação, Universidade de São Paulo, Avenida Trabalhador São-Carlense, 400, Centro, São Carlos, SP, 13566590, Brazil.
Dpto. Análisis Matemático y Matemática Aplicada, Universidad Complutense de Madrid, 28040, Madrid, Spain.

**Abstract**

We investigate some chemostat models incorporating wall growth, competition, random fluctuations on the dilution rate, and different consumption functions (Monod and Haldane). We analyze the asymptotic behavior of the solutions of the corresponding random differential systems to establish conditions on the model parameters under which the microbes persist in the gut or disappear. Moreover, several numerical simulations are presented to support the theoretical results and illustrate their biological interpretation.

## 1. Introduction

The human microbiome encompasses the diverse community of microorganisms, bacteria, viruses, fungi, and other microbes, that inhabit the human body. These microbial populations are essential to numerous physiological processes, including digestion, immune regulation, and defense against pathogenic organisms.

Each individual harbors a distinct microbiome, shaped by a range of factors such as diet, environment, and lifestyle. It is essential for maintaining overall health, and disruptions to the microbiome can be linked to various diseases and conditions.

---

[1]corresponding author: javier.lopez.delacruz@upm.es

We will focus, in particular, on the human intestinal microbiome, specifically the microbes that reside in the human gut, motivated by the strong correlation between the structure of the intestinal microbiome and the emergence of various diseases (see [9, 10, 15, 19, 21, 24, 28, 29, 31, 35]), which has sparked significant interest within the scientific community.

There are numerous references in the literature (see [20, 25, 39]) that attempt to understand the evolution of the human intestinal microbiome and seek to develop models that provide valuable insights, emphasizing that it is a complex problem.

Unlike most of papers of the existing literature, and building on [22], this article aims to address the problem of understanding the evolution of the human microbiome using mathematical models based on differential equations. Specifically, we will use the chemostat model, since we believe it is well-suited to describing the dynamics of the biological phenomenon under investigation (see [23] and references therein).

The chemostat device, introduced independently in the 1950s by Monod [30] and Novick and Szilard [32], is a well-established framework for describing the growth of microorganisms consuming a single nutrient in a controlled environment. This device serves as an idealized representation of natural microbial ecosystems where species (in our case) compete for a common resource.

The chemostat enjoys a wide range of practical applications, including ecology [3, 11, 16, 26, 40, 41, 42, 44], studies of the mammalian large intestine [18], wastewater treatment [27, 12], and genetic engineering challenges involving recombinant organisms [17]. Comprehensive discussions on applications of the chemostat can be found in [13, 37, 36, 38, 34, 43].

A chemostat includes three recipients: a feed reservoir, a culture chamber, and a waste collection vessel, all linked by pumps. The substrate is continuously introduced from the first tank into the second one, where the microbial growth. Simultaneously, an equal volume of culture is removed from the chamber to the collection vessel, ensuring that the culture volume does not vary over time.

The classical chemostat model is given by the following differential system (see, for instance, [23, 34])

$$s'(t) = D(s_{\text{in}} - s(t)) - c\mu(s(t))m(t),$$
$$m'(t) = -Dm(t) + g\mu(s(t))m(t),$$

where $s = s(t)$ and $m = m(t)$ denote, respectively, the concentrations of the nutrient and the microorganisms in the culture vessel at time $t$. The parameter $D > 0$ represents both the inflow rate of the nutrient and the outflow rate of the medium, ensuring volume balance. The constant $s_{\text{in}} > 0$ specifies the nutrient concentration in the inflowing medium, thus indicating the amount of nutrient available to the system. The parameter $c > 0$ stands for the maximal nutrient consumption rate, while $g \in (0, c]$ measures the effective microbial growth rate associated with nutrient consumption. Finally, $\mu : [0, +\infty) \to [0, +\infty)$ is the consume function, which characterizes how the microorganisms consume the nutrient.

Regarding the consumption function of the species, it is reasonable to assume that:

1. $\mu$ is continuous;
2. $\mu(0) = 0$;
3. $\mu(m) > 0$, $\forall m > 0$;
4. The parameters $c$ and $g$ can be adjusted, then we can assume without lost of generality that $\mu(m) \leq 1$ for all $m \geq 0$.

Because microorganisms may persist in the growth medium for extended periods and undergo natural mortality, we introduce a parameter $d > 0$ to represent their death rate coefficient. Furthermore, under conditions of limited dilution, the accumulated dead biomass can be broken down by bacterial activity, thereby contributing to nutrient recycling. To capture this effect, we let $r \in (0, 1)$ denote the proportion of dead biomass that undergoes recycling, i.e., the portion converted back into available nutrient.

After incorporating the aforementioned ingredients, we obtain

$$s'(t) = D(s_{\text{in}} - s(t)) - c\mu(s(t))m(t) + rdm(t),$$
$$m'(t) = -Dm(t) - dm(t) + g\mu(s(t))m(t),$$

where the state variables and parameters have been already introduced.

At this stage, it is worth noting that microorganisms can sometimes proliferate within the culture medium as well as on the wall of the container (see [4]). This phenomenon, known as wall growth, was introduced by Pilyugin and Waltman (see [33]) and is frequently observed in practice, especially under conditions of low flow rate.

Thus, we will divide the population of microbes into two: the population floating in the growth medium, that we will denote by $m_1$, and the population on the walls of the container, that we will denote by $m_2(t)$. As expected, those species on the wall can switch from the wall to the liquid media to float, and the floating individuals can get attached to the wall as well.

Assuming that the nutrient is equally available for both categories and they consume the same amount of nutrient, we obtain the new chemostat model

$$s'(t) = D(s_{\text{in}} - s(t)) - c\mu(s(t))(m_1(t) + m_2(t)) + rdm_1(t), \tag{1}$$
$$m_1'(t) = -(D + d)m_1(t) + g\mu(s(t))m_1(t) - \alpha_1 m_1(t) + \alpha_2 m_2(t), \tag{2}$$
$$m_2'(t) = -dm_2(t) + g\mu(s(t))m_2(t) + \alpha_1 m_1(t) - \alpha_2 m_2(t), \tag{3}$$

where $\alpha_1 \geq 0$ and $\alpha_2 \geq 0$ are the rates of microbial attachment to and detachment from the walls, respectively.

Next, we will introduce competition terms between the species, motivated by [22]. This situation usually arises due to the limitation of nutrients and need to be taken into account since our goal is to obtain models as much realistic as possible. Moreover, the classical chemostat model assumes that the dilution rate $D$ is the same for both entry and exit, which makes sense when considering a sufficiently large culture vessel (such as in industrial applications). However, this is not the case here, so we need to introduce a new parameter $\alpha$. Then, the resulting chemostat models is given as

$$s'(t) = D(s_{\text{in}} - \alpha s(t)) - c\mu(s(t))(m_1(t) + m_2(t)) + rdm_1(t), \tag{4}$$
$$m'_1(t) = m_1(t)\left(-d - \alpha D + g\mu(s(t)) - r_1 y_1(t) - r_2 m_2(t) - \alpha_1\right) + \alpha_2 m_2(t), \tag{5}$$
$$m'_2(t) = m_2(t)\left(-d + g\mu(s(t)) - r_1 m_1(t) - r_2 m_2(t) - \alpha_2\right) + \alpha_1 m_1(t), \tag{6}$$

where $\alpha > 0$ is the ratio between output and input flow rates, $r_1$ and $r_2$ denote the intensity of competition among microbial populations within the liquid medium and on the wall of the gut, respectively.

Although the competition terms makes the model much more realistic, some assumptions must still be made, despite being rather stringent from a biological standpoint. For example, we assume that the input flow $D$ is constant, despite the fact that it is often subject to bounded random fluctuations in real devices (see, for instance [6, 8]).

Although incorporating competition terms renders the model more realistic, it remains necessary to impose certain assumptions that are rather strong from a biological perspective. For example, we assume a constant input flow $D$, even though in actual devices this parameter is often exposed to bounded random fluctuations (see, e.g., [6, 8]).

This idea generalizes the work made in [22], where the authors studied the case where $D$ is not constant but deterministic. More precisely, we not only study a much more realistic system, but additionally we are able to improve some of their results.

Building on these considerations, we aim to conduct a thorough analysis of the long-term dynamics of the system (4)-(6), taking into account the influence of random fluctuations in the input and output flow. Nevertheless, it is essential to emphasize that the perturbed input/output flow is bounded in real devices, see e.g. [2, 14]. Hence, we will model the random fluctuations employing a bounded noise derived from the well-known Ornstein-Uhlenbeck process, widely recognized as a reliable approach for representing bounded stochastic disturbances (see, for instance, [8]).

This leads to the following formulation of the random chemostat model

$$s'(t) = (D + \psi(\xi^*(\theta_t\omega)))(s_{\text{in}} - \alpha s(t)) - c\mu(s(t))(m_1(t) + m_2(t)) + rdm_1(t), \tag{7}$$
$$m'_1(t) = m_1(t)\left(-d - \alpha(D + \psi(\xi^*(\theta_t\omega))) + g\mu(s(t)) - r_1 m_1(t) - r_2 m_2(t) - \alpha_1\right)$$
$$\qquad + \alpha_2 m_2(t), \tag{8}$$
$$m'_2(t) = m_2(t)\left(-d + g\mu(s(t)) - r_1 m_1(t) - r_2 m_2(t) - \alpha_2\right) + \alpha_1 m_1(t), \tag{9}$$

where $\xi^*$ is the Ornstein-Uhlenbeck process, see (11), and $\psi$ is an auxiliary function designed to ensure that the noise remains bounded, see (12).

Thus, our main goal is to find conditions on the model parameters under which the microbes persist in the gut or disappear. In addition, to make this manuscript more comprehensive, we will consider two consumption functions of the species. On the one hand, the Michaelis-Menten (also called Monod) consumption function

$$\mu_1(m) = \frac{m}{k+m}, \tag{10}$$

where $k > 0$ represents the half saturation constant (see [7, 34]), and, on the other hand, the Haldane consumption function

$$\mu_2(m) = \frac{m}{k+m+\frac{m^2}{i}},$$

which is non-monotonic, where $i > 0$ denotes the inhibition constant (see [7]).

Since real experiments involve bounded perturbations, the present work contributes to strengthening the biological interpretation of the ideal (unperturbed) system (4)–(6), as we will prove that it exhibits similar behavior under small perturbations. In other words, under conditions analogous to those of the non-perturbed model, extinction and persistence results for the species can still be obtained even in the presence of bounded perturbations. This justifies the suitability of the model for use in real experimental settings.

The subsequent sections of the paper are arranged as follows. In Section 2 we include some preliminaries on the use of the Ornstein-Uhlenbeck process to represent bounded stochastic fluctuations. Section 3 is devoted to proving the existence and uniqueness of a non-negative global solution to the random chemostat system. (7)-(9). Furthermore, we derive results regarding the existence of sets that absorbs and attract solutions. After that, in Section 4 we analyze in detail the long-time dynamics of the solutions of system (7)-(9). More precisely, we prove some results providing conditions on the parameters of the corresponding systems under which the microbes persist in the gut or disappear. Section 5 provides numerical simulations to corroborate the theoretical findings and highlight their biological significance. In Section 6 we conclude the paper by highlighting the main results and their scientific relevance.

## 2. Preliminaries about the bounded noise

We include below some preliminaries about the bounded noise $\xi^*(z^*(\theta_t\omega))$, necessary to facilitate the understanding of the rest of the paper and make it as much self-contained as possible.

Let us start considering the probability space $(\Omega, \mathcal{F}, \mathbb{P})$, where $\Omega$ is the space of every real continuous function being zero at zero, $\mathcal{F}$ is the Borel $\sigma-$algebra on $\Omega$ (see Appendix A.2 and Appendix A.3 in [1]) and $\mathbb{P}$ the corresponding Wiener measure.

In addition, let consider $\{\theta_t : \Omega \to \Omega\}_{t \in \mathbb{R}}$, where

$$\theta_t \omega(\cdot) = \omega(\cdot + t) - \omega(t),$$

which is known as Wiener shift flow.

Next, define the following stochastic process

$$\xi(t, \omega) := \xi^*(\theta_t \omega) = -\int_{-\infty}^{0} e^s \theta_t \omega(s) ds, \quad t \in \mathbb{R}, \ \omega \in \Omega, \tag{11}$$

which is mean-reverting (i.e., the probability of the process to go back to its mean value increases when the process is far away from its mean value) with continuous trajectories. See [5] for more details and properties.

Then, we can consider now the mapping $\psi : \mathbb{R} \to [-a, a]$, given by

$$\psi(\xi) = \frac{2a}{\pi} \arctan(\xi) \tag{12}$$

where $a > 0$ is a constant typically provided by practitioners, and define a new stochastic process as $\psi(\xi^*(\theta_t \omega))$, which satisfies some useful properties enumerated below.

**Proposition 2.1.** *Let $\psi$ be a function given as in (12) and consider $\xi^*(\theta_t \omega)$ the Ornstein-Uhlenbeck process. Then:*

1. *$\psi(\xi^*(\theta_t \omega))$ has continuous trajectories for almost every $\omega \in \Omega$,*

2. *$\psi(\xi^*(\theta_t \omega))$ is bounded for every $\omega \in \Omega$, in fact,*

$$D_{min} \leq D + \psi(\xi^*(\theta_t \omega)) \leq D_{max},$$

   *where $D_{min} := D - a$ and $D_{max} = D + a$, and*

3. *the following property*

$$\lim_{t \to +\infty} \frac{1}{t} \int_0^t \psi(\xi^*(\theta_s \omega)) ds = 0, \quad a.s. \ in \ \Omega. \tag{13}$$

   *fulfills for almost every $\omega \in \Omega$.*

### 3. Existence and uniqueness of solution. Absorbing and attracting sets

Henceforth, we denote $\mathbb{R}^3_+ := \{(x, y, z) \in \mathbb{R}^3 : x, y, z \geq 0\}$.

**Theorem 3.1.** *For any initial condition $u_0 = (s_0, m_{10}, m_{20}) \in \mathbb{R}^3_+$, system (7)-(9) has a unique global solution which lies in $\mathbb{R}^3_+$.*

*Proof.* Since the terms on the right-hand side of (7)-(9) are continuous with respect to time $t$ and $C^1$ (hence locally Lipschitz) as a function of $(s, m_1, m_2)$, we already have the existence and uniqueness of local solutions. That is, there exists a maximal time $t_M > 0$ such that $s, m_1, m_2$ are defined on $[0, t_M)$, satisfying the initial condition. To ensure that the solution lies in $\mathbb{R}_+^3$, note that if there exists $t_0 \in (0, t_M)$ such that $s(t_0) = 0$, then from (7), we have

$$0 < (D + \psi(\xi^*(\theta_t\omega)))s_{\text{in}} + rdm_1(t) = s'(t_0).$$

Hence, for each $t \in [0, t_M)$ it holds $s(t) \geq 0$. Analogously, if $m_1(t_0) = 0$ then

$$0 \leq m_1'(t_0) = \alpha_2 m_2$$

and $m_1(t) \geq 0$ for all $t \in [0, t_M)$. Doing the same for $m_2$ we obtain that

$$(s(t), m_1(t), m_2(t)) \in \mathbb{R}_+^3, \forall t \in [0, t_M).$$

Now we will show that the solutions do not blow up in finite time, that is, they are globally defined. Indeed, defining

$$z = gs + c(m_1 + m_2) \tag{14}$$

it is easy to see that

$$\begin{aligned}z'(t) &= g(D + \psi(\xi^*(\theta_t\omega)))s_{\text{in}} + grdm_1(t) - (g\alpha(D + \psi(\xi^*(\theta_t\omega))))s \\ &\quad + cdm_1(t) + c\alpha(D + \psi(\xi^*(\theta_t\omega)))m_1(t)) \\ &\quad - (cdm_2(t) + cr_1 m_1^2(t) + cr_2 m_2^2(t) + cr_2 m_1(t)m_2(t) + cr_1 m_1(t)m_2(t)) \\ &\leq g(D + \psi(\xi^*(\theta_t\omega)))s_{\text{in}} + grdz.\end{aligned}$$

It follows from the the Gronwall inequality that $z$ remains bounded in finite time, which concludes the proof. □

In the next theorem we show the existence of a set that aborbs solutions. Notice that we do not demand any extra condition (see [22]).

**Theorem 3.2.** *System (7)-(9) has an absorbing set that does not depend on $\omega$. That is, there exists $R > 0$ such that if $B \subset \mathbb{R}_+^3$ is bounded, then there exists $t_B > 0$ satisfying*

$$\| z(t; 0, \omega, z_0) \| \leq M, \ \forall t \geq t_B, \ \forall z_0 \in B, \ \forall \omega \in \Omega,$$

*where $z$ is given in* (14).

*Proof.* Using the same notation from Theorem 3.1, we have

$$\begin{aligned}
z'(t) &\leq g(D + \psi(\xi^*(\theta_t\omega)))s_{\text{in}} - g(D + \psi(\xi^*(\theta_t\omega)))\alpha s(t) + grdm_1(t) - cdm_1(t) \\
&\quad - c\alpha(D + \psi(\xi^*(\theta_t\omega)))m_1(t) - cdm_2(t) \\
&\leq g(D + \psi(\xi^*(\theta_t\omega)))s_{\text{in}} - gD_{\min}\alpha s(t) + grdm_1(t) - cdm_1(t) - c\alpha D_{\min}m_1(t) - cdm_2(t) \\
&= g(D + \psi(\xi^*(\theta_t\omega)))s_{\text{in}} - gD_{\min}\alpha s(t) - (cd + c\alpha D_{\min} - grd)m_1(t) - cdm_2(t) \\
&\leq g(D + \psi(\xi^*(\theta_t\omega)))s_{\text{in}} - \vartheta z,
\end{aligned}$$

where
$$\vartheta = \min\{D_{\min}\alpha, d + \alpha D_{\min} - gc^{-1}rd, d\} > 0.$$

Therefore, for each $\omega \in \Omega$ and $z_0 \in \mathbb{R}^3_+$ it holds

$$z(t;0,\omega,z_0) \leq z_0 e^{-\vartheta t} + s_{\text{in}}g\int_0^t (D + \psi(\xi^*(\theta_t\omega)))e^{-\vartheta(t-s)}ds, \quad \forall t \geq 0$$

and consequently

$$z(t;0,\omega,z_0) \leq z_0 e^{-\vartheta t} + \frac{gD_{\max}s_{\text{in}}}{\vartheta}(1 - e^{-\vartheta t}), \quad \forall t \geq 0, \ \forall \omega \in \Omega, \forall z_0 \in \mathbb{R}^3_+.$$

Thus, given $\epsilon > 0$, the compact subset

$$\mathcal{B}_\epsilon := \left\{(s, m_1, m_2) \in \mathbb{R}^3_+ : 0 \leq gs + c(m_1 + m_2) \leq \frac{gD_{\max}s_{\text{in}}}{\vartheta} + \epsilon\right\} \tag{15}$$

is absorbing for system (7)-(9). □

**Remark 3.1.** *From Theorem 3.2, we proved that the set*

$$\mathcal{B}_0 := \left\{(s, m_1, m_2) \in \mathbb{R}^3+ : 0 \leq gs + c(m_1 + m_2) \leq \frac{gD_{max}s_{in}}{\vartheta}\right\}$$

*is a deterministic compact attracting set for the solutions of (7)-(9). Moreover, if $s_{in} \leq 1$, then $z'\left(\frac{gD_{max}}{\vartheta}\right) \leq 0$ and $\mathcal{B}_0$ is positively invariant.*

## 4. Asymptotic Behavior

In order to study the asymptotic behavior of the solutions we will define the following variables:
$$m = m_1 + m_2 \quad \text{and} \quad p = \frac{m_1}{m}.$$

With this new variables we obtain the following system:

$$s'(t) = (D + \psi(\xi^*(\theta_t\omega)))(s_{\text{in}} - \alpha s(t)) - c\mu(s(t))m(t) + rdp(t)m(t) \tag{16}$$

$$m'(t) = -\alpha(D + \psi(\xi^*(\theta_t\omega)))p(t)m(t) + g\mu(s(t))m(t) - dm(t) - r_1 m^2(t) - r_2(1-p(t))m^2(t) \tag{17}$$

$$p'(t) = \alpha(D + \psi(\xi^*(\theta_t\omega)))p^2(t) - (\alpha(D + \psi(\xi^*(\theta_t\omega))) + \alpha_1 + \alpha_2)p(t) + \alpha_2. \tag{18}$$

Fix $\epsilon > 0$ arbitrary. Since $p(t) \in [0, 1]$ for each $t \geq 0$, then

$$p'(t) \leq \alpha_2 - (\alpha_1 + \alpha_2)p(t)$$

and therefore

$$p(t) \leq p(0)e^{-(\alpha_1+\alpha_2)t} + \frac{\alpha_2}{\alpha_1 + \alpha_2}(1 - e^{-(\alpha_1+\alpha_2)t}).$$

Hence, we can fix $t_\epsilon > 0$ large satisfying

$$p(t) \leq P_\epsilon := \frac{\alpha_2}{\alpha_1 + \alpha_2} + \epsilon, \ \forall t \geq t_\epsilon.$$

Analogously, we have

$$p(t) \geq p(0)e^{-(\alpha D + \alpha_1 + \alpha_2)t - \alpha \int_0^t \psi(\xi^*(\theta_s\omega))ds} + \alpha_2 \int_0^t e^{-(\alpha D + \alpha_1 + \alpha_2)(t-s)}ds.$$

From the ergodic property (13) we can increase $t_\epsilon > 0$, if necessary, to ensure that

$$p(t) \geq p_\epsilon := \frac{\alpha_2}{\alpha_1 + \alpha_2 + \alpha D_{\max}} - \epsilon, \ \forall t \geq t_\epsilon.$$

This shows that for each $\epsilon > 0$ there exists $t_\epsilon > 0$ satisfying

$$p(t) \in \left(\frac{\alpha_2}{\alpha_1 + \alpha_2 + \alpha D_{\max}} - \epsilon, \frac{\alpha_2}{\alpha_1 + \alpha_2} + \epsilon\right), \ \forall t \geq t_\epsilon \tag{19}$$

and

$$p_\epsilon m_2(t) < m_1(t)(1 - p_\epsilon) < m_2(t)P_\epsilon \left(\frac{1 - p_\epsilon}{1 - P_\epsilon}\right), \ \forall t \geq t_\epsilon. \tag{20}$$

**Theorem 4.1** (Extinction). *Assume that*

$$D\alpha \frac{\alpha_2}{\alpha_1 + \alpha_2 + \alpha D_{max}} + d > g. \tag{21}$$

*Then the species will be extinct, meaning that the solutions of equation (17) will converge to zero as t approaches infinity.*

*Proof.* Using the same notation as on Section 3 fix $\epsilon > 0$ such that
$$D\alpha p_\epsilon + d > g.$$

From equation (13), (18) and (19), we obtain that for large $t$, the following inequality holds
$$m'(t) \leq \left(-\alpha(D + \psi(\xi^*(\theta_t\omega)))p_\epsilon + g - d\right)m(t),$$

which implies
$$m(t) \leq m(0)\exp\left\{-(\alpha D p_\epsilon - g + d)t - \alpha p_\epsilon \int_0^t \psi(\xi^*(\theta_\tau\omega))d\tau\right\}.$$

Therefore, $m(t) \to 0$ as $t \to +\infty$. □

Note that our extinction condition (21) is weaker than
$$D_{\min}\alpha \frac{\alpha_2}{\alpha_1 + \alpha_2 + \alpha D_{\max}} + d > g,$$

which is the condition obtained in [22]. This improvement further justifies our choice of bounded noise.

From this point onward, we aim to find conditions that guarantee the survival of the species. To achieve this, we will first determine a lower bound for the consumption function $s$.

**Proposition 4.1.** *There exist $T > 0$ and $s^* > 0$ satisfying*
$$s(t) \geq s^*, \ \forall t \geq T.$$

*Proof.* Fix $\epsilon > 0$ arbitrary. Applying Theorem 3.2 we obtain
$$cm \leq \frac{gD_{\max}s_{\text{in}}}{\vartheta} + \epsilon,$$

for large $t$. Hence
$$s'(t) \geq D_{\min}s_{\text{in}} - \alpha D_{\max}s(t) - \mu(s(t))\left(\frac{gD_{\max}}{\vartheta} + \epsilon\right), \quad \text{for large } t.$$

Define $f : [0, +\infty) \to [0, +\infty)$ as
$$f(s) = D_{\min}s_{\text{in}} - \alpha D_{\max}s - \mu(s)\left(\frac{gD_{\max}}{\vartheta} + \epsilon\right), \ \forall s \geq 0.$$

Since $f$ is continuous and $f(0) > 0$, we may choose $s^* > 0$ in such a way that $f(s) > 0$ for $s \in [0, s^*]$. Therefore,
$$s(t) > s^*, \quad \text{for large } t. \tag{22}$$
□

**Remark 4.1.** *Note that $s^*$ in Proposition 4.1 can be computed numerically, as it corresponds to the smallest positive root of the function $f$, as shown in the proof.*

### 4.1. Monod consumption function

In what follows, we derive conditions that ensures the survival of the species, taking into account the Monod consumption function given by (10).

**Theorem 4.2.** *Let $k > 0$ be the half-saturation constant given in (10) and define*

$$l = \min\{s^*, D_{max}k^{-1}\} > 0.$$

*If*

$$\alpha D_{max}\frac{\alpha_2}{\alpha_1 + \alpha_2} + d < g\frac{s^*}{k+l}, \tag{23}$$

*then the survival of the species is ensured, that is, there exists $m_1^*, m_2^*, t_s > 0$ such that*

$$m_1(t) > m_1^* \quad \text{and} \quad m_2(t) > m_2^* \; \forall t > t_s.$$

*Proof.* First we will prove that there exists $m^* > 0$ such that

$$m(t) > m^*, \text{ for large } t.$$

Let us consider the case $s^* \leq D_{\max}k^{-1}$. Fixing $\epsilon > 0$ such that

$$\alpha D_{\max} P_\epsilon + d < g\frac{s^*}{k+s^*}$$

it follows from (17), (22) and the fact that the Monod consumption function is increasing that

$$m'(t) \geq m(t)g(m(t)),$$

where $g : [0, +\infty) \to [0, +\infty)$ is the continuous function

$$g(m) = -\alpha D_{\max} P_\epsilon + g\frac{s^*}{k+s^*} - d - (r_1 P_\epsilon + r_2 p_\epsilon - r_2)m, \; \forall m \in [0, +\infty).$$

□

Since $g(0) > 0$ there exists $m^* > 0$ satisfying $g(t) > 0$ for all $t \in [0, m^*]$ and consequently

$$y(t) > m^* \quad \text{for large } t. \tag{24}$$

If $D_{\max}k^{-1} < s^*$ fix $\epsilon > 0$ such that

$$\alpha D_{\max} P_\epsilon + d < g\frac{s^*}{k + D_{\max}k^{-1} + \epsilon g^{-1}}.$$

From (15) we know that
$$s \leq D_{\max}k^{-1} + g^{-1}\epsilon - cg^{-1}m, \quad \text{for large } t$$

and therefore
$$m'(t) \geq m(t)h(m(t)), \quad \text{for large } t$$

where $h : [0, H] \to [0, +\infty)$ is the continuous function

$$h(m) = -\alpha D_{\max}P_\epsilon + g\frac{s^*}{k + D_{\max}k^{-1} + g^{-1}\epsilon - cg^{-1}m} - d - (r_1 P_\epsilon + r_2 p_\epsilon - r_2)m, \quad \forall y \in [0, H]$$

and $H = \frac{g}{c}(k + D_{\max}k^{-1} + g^{-1}\epsilon)$. Since $h(0) = 0$ we can argue as the another case to obtain (24).

It follows from (19) and (20) that

$$\liminf_{t \to +\infty} m_1(t) > \left(\frac{\alpha_2}{\alpha_1 + \alpha_2 + \alpha D_{\max}}\right)m^* \tag{25}$$

and

$$\liminf_{t \to +\infty} m_2(t) > \left(\frac{\alpha_2}{\alpha_1 + \alpha_2 + \alpha D_{\max}}\right)\left(\frac{\alpha_2}{\alpha_1 + \alpha_2}\right)m^*. \tag{26}$$

**4.2. Haldane consumption function**

In this subsection, we will examine the Haldane consumption function, defined as

$$\mu(x) = \frac{s}{k + s + \frac{s^2}{i}}, \quad \forall s \geq 0,$$

where $i > 0$ is the inhibition constant [7]. The following result guarantees the persistence of the species under the Haldane consumption function.

**Theorem 4.3.** *Assume that*

$$\alpha D_{max}\frac{\alpha_2}{\alpha_1 + \alpha_2} + d < \frac{s^*}{k + k^{-1}D_{max}s_{in} + D_{max}^2 s_{in}^2 i^{-1}k^{-2}}. \tag{27}$$

*Then the survival of the species is ensured, that is, there exists $m_1^*, m_2^*, t_s > 0$ such that*

$$m_1(t) > m_1^* \quad \text{and} \quad m_2(t) > m_2^* \ \forall t > t_s.$$

*Proof.* Since the proof follows similarly to that of Theorem 4.2, we will present it with fewer details. Fix $\epsilon > 0$ such that

$$\alpha D_{\max}\frac{\alpha_2}{\alpha_1 + \alpha_2} + d < \frac{s^*}{k + (k^{-1}D_{\max}s_{\text{in}} + \epsilon) + i^{-1}(D_{\max}s_{\text{in}}k^{-1} + \epsilon)^2}.$$

From (15) we have
$$s(t) < \frac{D_{\max} s_{\text{in}}}{k} + \epsilon, \quad \text{for large } t.$$

Hence
$$m'(t) \geq m(t)(-\alpha D_{\max} P_\epsilon - d + g\mu(s(t))) + m^2(t)(p_\epsilon - r_1 - r_2)$$
$$\geq m(t)c(m(t)),$$

where
$$c(m) = \left(-\alpha D_{\max} P_\epsilon - d + g\frac{s^*}{k + (k^{-1}D_{\max}s_{\text{in}} + \epsilon) + i^{-1}(D_{\max}s_{\text{in}}k^{-1} + \epsilon)^2}\right) + m(p_\epsilon - r_1 - r_2)\right).$$

Since $c$ is continuous and $c(0) > 0$, we can proceed as in Theorem 4.2 to guarantee the existence of a value $m_h^* > 0$ satisfying
$$m(t) > m_h^*, \quad \text{for large } t$$

and consequently
$$\liminf_{t \to +\infty} m_1(t) > \left(\frac{\alpha_2}{\alpha_1 + \alpha_2 + \alpha D_{\max}}\right) m_h^* \tag{28}$$

and
$$\liminf_{t \to +\infty} m_2(t) > \left(\frac{\alpha_2}{\alpha_1 + \alpha_2 + \alpha D_{\max}}\right)\left(\frac{\alpha_2}{\alpha_1 + \alpha_2}\right) m_h^*. \tag{29}$$

□

## 5. Numerical simulations

Numerical simulations are carried out in this section to demonstrate and validate the analytical results discussed in Section 4.

Each figure consists of three panels. The top panel shows the time evolution of the nutrient concentration. The blue dashed line represents the solution of system (4)-(6), while the continuous colored lines correspond to the solutions of the random system (7)-(9) for different realizations of the noise. Similarly, the middle panel displays the time evolution of the species concentration in the medium, and the bottom panel shows the evolution of the species concentration on the wall.

In Figure 1 and Figure 2, we illustrate the cases with both consumption functions, respectively. For these simulations, we set the following parameters: $c = 4.8$, $g = 0.6$, $k = 4.7$, $s_{\text{in}} = 17$, $\alpha = 0.5$, $r = 0.4$, $d = 0.4$, $\alpha_1 = 0.5$, $\alpha_2 = 0.7$, $D = 1.9$, $a = 0.25$, $r_1 = 0.1$, $r_2 = 0.5$ and the initial condition $(s_0, m_{10}, m_{20}) = (20, 14, 10)$.

In both cases, it is straightforward to verify that condition (21) holds, leading to the extinction of the species (both those in the medium and on the wall), as shown in the previous figures. It should be emphasized that the dynamics of the random model (7)-(9) is similar for both the Monod and Haldane cases. Thus, the inhibition constant in the Haldane case does not significantly affect the dynamics when condition (21) holds.

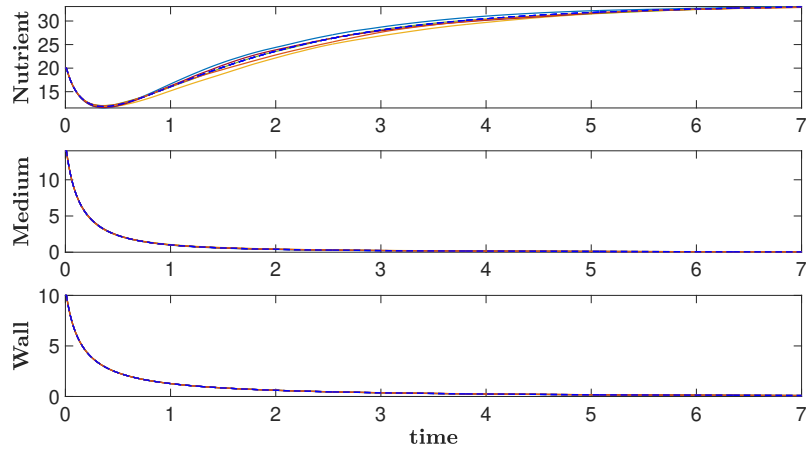

Figure 1: **Extinction (Monod case)**. Parameters: $c = 4.8$, $g = 0.6$, $k = 4.7$, $s_{\text{in}} = 17$, $\alpha = 0.5$, $r = 0.4$, $d = 0.4$, $\alpha_1 = 0.5$, $\alpha_2 = 0.7$, $D = 1.9$, $a = 0.25$, $r_1 = 0.1$, $r_2 = 0.5$. Initial condition: $(s_0, m_{10}, m_{20}) = (20, 14, 10)$.

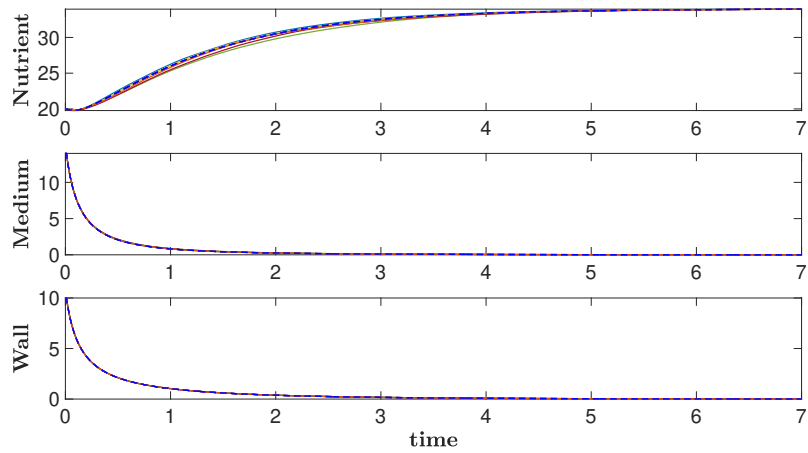

Figure 2: **Extinction (Haldane case)**. Parameters: $c = 4.8$, $g = 0.6$, $k = 4.7$, $s_{\text{in}} = 17$, $\alpha = 0.5$, $r = 0.4$, $d = 0.4$, $\alpha_1 = 0.5$, $\alpha_2 = 0.7$, $i = 5$, $D = 1.9$, $a = 0.25$, $r_1 = 0.1$, $r_2 = 0.5$. Initial condition: $(s_0, m_{10}, m_{20}) = (20, 14, 10)$.

Next, in Figure 3 (using the Monod consumption function), we consider the parameters $c = 4$, $g = 3.8$, $k = 1.4$, $s_{\text{in}} = 17$, $\alpha = 0.5$, $r = 0.4$, $d = 0.01$, $\alpha_1 = 0.5$, $\alpha_2 = 0.7$, $D = 0.47$, $a = 0.4$, $r_1 = 0.4$, $r_2 = 0.6$, and the initial condition $(s_0, m_{10}, m_{20}) = (20, 14, 10)$. In this case, condition (23) is satisfied, ensuring species survival (both those floating and on the wall), as stated in Theorem 4.2 as observed in the simulation.

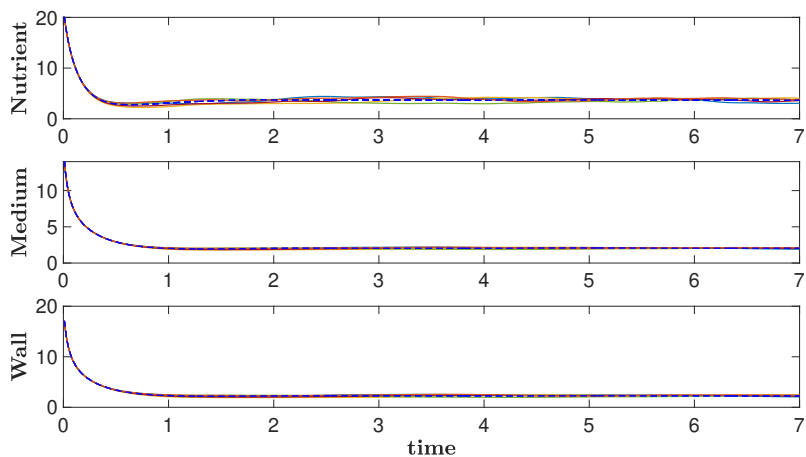

Figure 3: **Persistence (Monod case)**. Parameters: $c = 4$, $g = 3.8$, $k = 1.4$, $s_{\text{in}} = 17$, $\alpha = 0.5$, $r = 0.4$, $d = 0.01$, $\alpha_1 = 0.5$, $\alpha_2 = 0.7$, $D = 0.47$, $a = 0.4$, $r_1 = 0.4$, $r_2 = 0.6$. Initial condition: $(s_0, m_{10}, m_{20}) = (20, 14, 10)$.

Finally, in Figure 4 (using the Haldane consumption function), we consider the parameters $c = 7$, $g = 6.8$, $k = 7$, $s_{\text{in}} = 19$, $\alpha = 0.5$, $r = 0.7$, $d = 0.01$, $\alpha_1 = 0.55$, $\alpha_2 = 0.65$, $i = 7.6$, $D = 0.61$, $a = 0.16$, $r_1 = 0.4$, $r_2 = 0.2$, and the initial condition $(s_0, m_{10}, m_{20}) = (20, 14, 17)$. Since condition (27) is satisfied, the persistence of the species (both those in the medium and on the wall) is guaranteed, as demonstrated in Theorem 4.3, which is evident in the corresponding figure.

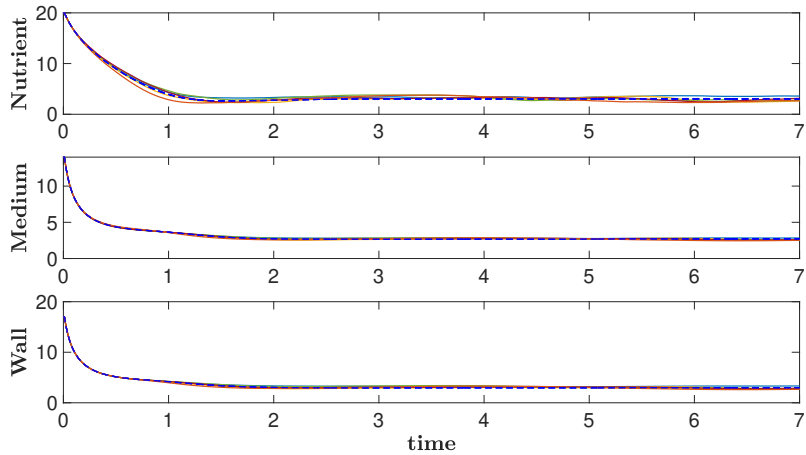

Figure 4: **Persistence (Haldane case)**. Parameters: $c = 7$, $g = 6.8$, $k = 7$, $s_{\text{in}} = 19$, $\alpha = 0.5$, $r = 0.7$, $d = 0.01$, $\alpha_1 = 0.55$, $\alpha_2 = 0.65$, $i = 7.6$, $D = 0.61$, $a = 0.16$, $r_1 = 0.4$, $r_2 = 0.2$. Initial condition: $(s_0, m_{10}, m_{20}) = (20, 14, 17)$.

## 6. Conclusions

In this work, we analyzed the long-time dynamics of microorganisms in the chemostat model with wall growth, noise, and competition, as described by system (7)-(9).

We derived conditions for both the extinction and survival of the microorganisms. In summary, we proved the following results:

1. If condition (21) is satisfied, the microorganisms will become extinct, regardless of the consumption function.

2. For the Monod consumption function, condition (23) guarantees the survival of the species, with lower bounds given by (25).

3. For the Haldane consumption function, condition (27) ensures the survival of the species, with lower bounds provided by (28).

Finally, we emphasize the strengths of our paper below:

- We rigorously prove that system (4)–(6) is stable under perturbations, in the following sense: even with bounded perturbations on the input/output flow $D$, the system still behaves similarly to the unperturbed one.

- The extinction condition (21) depends continuously on the perturbation magnitude $d$, so that when $d = 0$, it reduces to the extinction condition for the ideal system (4)–(6). This shows that if the extinction condition for the ideal (unperturbed) model (4)–(6) is satisfied, then extinction will also occur in the perturbed case, provided the perturbation is sufficiently small.

- Conditions (23) and (27), that guarantee species survival with the Monod and Haldane consumption functions respectively, depend continuously on $d$. Similar to the extinction case, if one of these conditions is satisfied for the non-perturbed case (where $d = 0$), then the survival of the species is still guaranteed in the perturbed case, provided that the perturbation is sufficiently small.

- We provide numerical simulations that corroborate all the theoretical results presented above and assist in visualizing the system's behavior.

- The above points justify that system (4)–(6) is suitable for modeling real-world problems, where bounded perturbations occur in the flow $D$.

**Financial disclosure**

This work has been partially supported by the Spanish Ministerio de Ciencia, Innovación y Universidades, Agencia Estatal de Investigación (AEI) and Fondo Europeo de Desarrollo Regional (FEDER) under the project PID2024-156228NB-I00, the Brazilian Coordination for the Improvement of Higher Level Personnel (CAPES/PROEX) under project D-11169228 in addition to the São Paulo Research Foundation (FAPESP) under projects 2020/14353-6, 2022/02172-2 and 2024/16879-6.